\documentclass{article}
\pdfoutput=1
\usepackage{amssymb,amsmath}
\usepackage{amsfonts}
\usepackage{latexsym}
\newtheorem{Lem}{Lemma}
\newtheorem{Thm}{Theorem}

\newtheorem{Rem}{Remark}
\newtheorem{Cor}{Corollary}
\newtheorem{Exa}{Example}
\newtheorem{Def}{Definition}
\def\C{{\mathbb C}}
\def\R{{\mathbb R}}

\def\H{{\mathbb H}}
\def\M{{\mathbb M}} 

\def\q{{\bf q}}
\def\s{{\bf s}}
\def\z{{\bf z}} 
\def\V{{\mathcal{V}}}
\def\A{{\mathcal{A}}}
\def\X{{\mathfrak{X}}}

\title{Spectrum and Analytic Functional Calculus in Real and Quaternionic Frameworks} 
\author{Florian-Horia Vasilescu\\
\small Department of Mathematics, University of Lille,\\
\small 59655 Villeneuve d'Ascq, France\\
\small florian.vasilescu@univ-lille.fr}
\date{} 

\begin{document}

\maketitle
\bigskip



\begin{abstract} 
We present an approach to the spectrum and analytic functional calculus 
for quaternionic linear operators, following the corresponding results 
concerning the real linear operators. In fact, the construction of the analytic functional calculus for real linear operators can be refined to get a similar construction for quaternionic linear ones, in a classical manner, using a Riesz-Dunford-Gelfand type kernel, and  considering spectra
in the complex plane. A quaternionic joint spectrum for
pairs of operators is also discussed, and an analytic functional calculus is constructed, via a Martinelli type kernel in two variables. 
\end{abstract}
\medskip

calculus. 


\section{Introduction}\label{I}

In this text we consider $\R$-,\,$\C$-, and $\H$-linear operators, that is, 
real, complex and quaternionic linear operators, respectively.

While the spectrum of a linear operator is traditionally defined 
for complex linear operators, it is sometimes useful to have it 
also for real linear operators, as well as for quaternionic linear ones. 
The definition of the spectrum
for a real linear operator goes seemingly back to Kaplansky (see \cite{Kap}), and it can be stated as follows. If $T$ is a real linear  operator on the real vector space $\V$, a point $u+iv$ 
($u,v\in\R$) is in the spectrum of $T$ if the operator $(u-T)^2+v^2$ is not 
invertible on $\V$, where the scalars are identified with multiples of the identity on $\V$. Although this definition involves only operators acting in $\V$, the spectrum is, nevertheless, a subset of the complex plane. As a matter of fact, a motivation of this choice can be illustrated via the complexification of the space $\V$ (see Section \ref{SC}). 

The spectral theory for quaternionic linear operators is largely
discussed in numerous work, in particular in the monographs     
\cite{CoSaSt} and \cite{CoGaKi}, and in many of their references as well.
In these works, the construction of an analytic functional
calculus (called $S$-{\it analytic functional calculus}) means to associate to each function from the class of the so-called {\it slice hyperholomorphic} or {\it slice regular functions} a quaternionic linear operator, using a specific noncommutative kernel.

The idea of the present work is to replace the class of slice regular functions by a class  holomorphic 
functions, using a commutative kernel of type Riesz-
Dunford-Gelfand. These two classes are isomorphic via a Cauchy type
transform (see \cite{Vas5}), and  the image of the analytic functional calculus is the same, as one might expect (see Remark \ref{twofc}).
 
As in the case of real operators, the verbatim extension of the classical definition of the spectrum for quaternionic operators is not appropriate, and so a different definition using the squares of operators and real numbers was given, which  can be found in \cite{CoSaSt} (see also \cite{CoGaKi}). We discuss this definition in our framework (see Definition \ref{Q-spectrum}), showing later that its ''complex border`` contains the most significant information, leading to the construction of an analytic functional calculus, equivalent to that obtained via the slice hyperholomorphic functions. 

In fact, we first consider the spectrum for real operators on real Banach spaces, and sketch the construction of an analytic
functional calculus for them, using some classical ideas (see Theorem \ref{R_afc}). Then we 
extend this framework to a quaternionic one, showing that the 
approach from the real case can be easily adapted to the new situation. 

As already mentioned, and unlike in  \cite{CoSaSt} or \cite{CoGaKi}, our functional calculus is obtained via a 
Riesz-Dunford-Gelfand formula, defined in a partially commutatative context, rather than the non-commutative Cauchy type formula used by previous authors. Our analytic functional calculus 
holds for a class of analytic operator valued functions, whose definition extends that of stem functions, and it applies, in particular, to a large family of quaternionic linear operators. Moreover, we can show that the analytic functional calculus obtained in this way is equivalent to the analytic functional calculus obtained in \cite{CoSaSt} or \cite{CoGaKi},
in the sense that the images of these functional calculi coincide (see Remark \ref{twofc}).

We finally discuss the case of pairs of commuting real operators,
in the spirit of \cite{Vas4}, showing some connections with the quaternionic case. Specifically, we define a quaternionic spectrum for them and construct an analytic functional calculus using a Martinelli type formula, showing that for such a construction only a sort of ''complex border`` of the 
quaternionic spectrum should be used. 

This work is just an introductory one. Hopefully, more contributions on this line will be presented in the future.

\section{Spectrum and Conjugation}\label{SC}

Let $\A$ be a unital real Banach algebra, not necessarily commutative. As mentioned in the Introduction,  the (complex) spectrum of an element $a\in\A$ may be defined by the equality

\begin{equation}\label{resp}
\sigma_\C(a)=\{u+iv;(u-a)^2+v^2\,\, {\rm is\,\,not\,\,invertible}, u,v\in\R\},
\end{equation} 

This set is {\it conjugate symmetric}, that is $u+iv\in\sigma_\C(a)$ if and only if $u-iv\in\sigma_\C(a)$. A known motivation of this definition comes from the following remark.

Fixing a unital real Banach algebra $\A$,  we denote by
$\A_\C$ the complexification of $\A$, which is given by
 $A_\C=\C\otimes_\R\A$, written simply as $\A+i\A$, where the sum is direct, identifying the element $1\otimes a+i\otimes b$ with the element  $a+ib$, for all $a,b\in\A$. 

 Then $\A_\C$ is a unital complex
algebra, which can be organized as a Banach algebra, with a (not
necessarily unique) convenient norm. To fix the ideas, we recall
that the product of two elements  is given by
$(a+ib)(c+id)=ac-bd+i(ad+bc)$ for all $a,b,c,d\in\A$, and the norm may be defind by
$\Vert a+ib\Vert=\Vert a\Vert+\Vert b\Vert$, where $\Vert*\Vert$
is the norm of $\A$.

In the algebra $\A_\C$,
the complex numbers commute with all elements of $\A$. Moreover,
we have a {\it conjugation} given by
$$
\A_\C\ni a+ib\mapsto a-ib\in\A_\C,\,a,b\in\A,
$$
which is a unital conjugate-linear automorphism, whose square is the identity. In particular, an arbitrary element $a+ib$ is invertible if and only if $a-ib$ is invertible. 

The usual spectrum, defined for each element $a\in\A_\C$, will be
denoted by $\sigma(a)$. Regarding the algebra $\A$ as a real
subalgebra of $\A_\C$, one has the following.

% Lem 1

\begin{Lem}\label{com_sp_re_op} For every $a\in\A$ we have the equality 
$\sigma_\C(a)=\sigma(a)$.
\end{Lem}

{\it Proof.} The result is well known but we give a short proof, because
a similar idea will be later used.

Let $\lambda=u+iv$ with $u,v\in\R$ arbitrary. Assuming 
$\lambda-a$  invertible, we also have $\bar{\lambda}-a$  invertible. From the obvious identity
$$
(u-a)^2+v^2=(u+iv-a)(u-iv-a),
$$
we deduce that the element $(u-a)^2+v^2$ is invertible, 
implying  the inclusion $\sigma_\C(a)\subset\sigma(a)$. 

Conversely, if $(u-a)^2+v^2$ is invertible, then both
$u+iv-a,u-iv-a$ are invertible via the decomposition  from above, showing that we also have $\sigma_\C(a)\supset\sigma(a)$.

% Rem 1

\begin{Rem}\rm The spectrum $\sigma(a)$ with $a\in\A$ is always a conjugate symmetric set.
\end{Rem}
 
We are particularly interested to apply the discussion from above to the context of linear operators. The spectral theory for real linear operators is well known, and it is developed actually  in the framework of linear relations (see \cite{BaZa}).
Nevertheless, we present here a different approach, which can 
be applied, with minor changes, to the case of some quaternionic operators.  

For a real or complex Banach space $\mathcal{V}$, we denote by
$\mathcal{B(V)}$ the algebra of all bounded $\R$-(
respectively $\C$-)linear operators on $\mathcal{V}$.
As before, the multiples of the identity will be identified with the corresponding scalars. 
 
Let $\V$ be a real Banach space, and let $\V_\C$ be its complexification, which, as above, is identified with the direct sum $\V+i\V$. Each operator
$T\in\mathcal{B(V)}$ has a natural extension to an operator
$T_\C\in\mathcal{B}(\V_\C)$, given by $T_\C(x+iy)=Tx+iTy,\,
x,y\in\V$. Moreover, the map $\mathcal{B(V)}\ni T\mapsto T_\C\in\mathcal{B}(\V_\C)$ is  unital, $\R$-linear and multiplicative. In particular, $T\in\mathcal{B(V)}$ is invertible if and only if
$T_\C\in\mathcal{B}(\V_\C)$ is invertible. 

Fixing an operator $S\in\mathcal{B}(\mathcal{V}_\C)$, we define 
the operator $S^\flat\in\mathcal{B}(\mathcal{V}_\C)$ to be equal to $CSC$, where $C:\V_\C\mapsto\V_\C$ is the conjugation
$x+iy\mapsto x-iy,\,x,y\in\V$.  It is easily seen that the map 
$\mathcal{B}(\mathcal{V}_\C)\ni S\mapsto S^\flat\in \mathcal{B}(\mathcal{V}_\C)$ is a unital conjugate-linear automorphism, whose square is the identity on $\mathcal{B}(\mathcal{V}_\C)$. 
Because $\mathcal{V}=\{u\in\mathcal{V}_\C; Cu=u\}$, we have
$S^\flat=S$ if and only if $S(\mathcal{V})\subset\mathcal{V}$.
In particular,  we have $T_\C^\flat=T_\C$. In fact, because of the representation
$$
S=\frac{1}{2}(S+S^{\flat})+ i\frac{1}{2i}(S-S^{\flat}),\,\,
S\in\mathcal{B}(\mathcal{V}_\C),
$$ 
where $(S+S^{\flat})(\V)\subset\V, i(S-S^{\flat})(\V)\subset\V$,  the algebras $\mathcal{B}(\mathcal{V}_\C)$ and $\mathcal{B(V)}_\C$ are isomorphic and they will be often identified, and $\mathcal{B(V)}$ will be regarded as a (real) subalgebra of $\mathcal{B}(\mathcal{V})_\C$. In particular, if $S=U+iV$, with $U,V\in\mathcal{B(V)}$, we have
$S^\flat=U-iV$, so the map $S\mapsto S^\flat$ is the  conjugation
of the complex algebra $\mathcal{B}(\mathcal{V})_\C$ induced
by the conjugation $C$ of $\V_\C$. 
 
For every operator $S\in\mathcal{B}(\mathcal{V}_\C)$, we denote,
as before, by $\sigma(S)$  its usual spectrum. As 
$\mathcal{B(V)}$ is a real algebra,
the (complex) spectrum of an operator $T\in\mathcal{B(V)}$ is given by the equality (\ref{resp}):
$$
\sigma_\C(T)=\{u+iv;(u-T)^2+v^2\,\, {\rm is\,\,not\,\,invertible}, u,v\in\R\}.
$$ 

\begin{Cor} For every $T\in\mathcal{B}(\V)$ we have the equality 
$\sigma_\C(T)=\sigma(T_\C)$.
\end{Cor} 

\section{Analytic Functional Calculus for Real 
\\ Operators}

Having a concept of spectrum for real operators, an important step for further development is the construction of an analytic functional calculus. Such a construction has been done actually
in the context of real linear relations in \cite{BaZa}. In
what follows we shall present a similar construction for real
linear operators. Although the case of linear relations looks more general,  unlike in \cite{BaZa}, we perform our construction using a class of operator valued analytic functions insted of scalar valued  analytic functions. Moreover, our arguments look simpler, and 
the construction is a model for a more general one, to get an analytic 
functional calculus for quaternionic linear operators.

If $\V$ is a real Banach space, and so each operator $T\in\mathcal{B}(\V)$ has a complex spectrum $\sigma_\C(T)$, which is
compact and nonempty, one can  use the classical Riesz-Dunford functional calculus, in a slightly generalized form (that is, replacing the scalar-valued analytic functions by operator-valued analytic ones, which is a well known idea).

The use of  vector versions of the Cauchy formula is simplified by adopting the following definition. Let $U\subset\C$ be open. An open subset 
$\Delta\subset U$ will be called a {\it Cauchy domain} (in $U$) if 
$\Delta\subset\bar{\Delta}\subset U$ and the boundary of $\Delta$ consists of a finite family of closed curves, piecewise smooth, positively oriented. Note that a Cauchy domain is bounded but not necessarily connected.

% Rem 2

\begin{Rem}\label{afcro}\rm If $\mathcal{V}$ is a real Banach space, and 
$T\in\mathcal{B(V})$, we have the usual analytic functional calculus for the operator $T_\C\in\mathcal{B}(\mathcal{V}_\C)$
(see \cite{DuSc}). That is, in a slightly generalized form, and for later use, if $U\supset\sigma(T_\C)$ is an open set in $\C$ and 
$F:U\mapsto\mathcal{B}(\mathcal{V}_\C)$ is analytic, we  put
$$
F(T_\C)=\frac{1}{2\pi i}\int_\Gamma F(\zeta)(\zeta-T_\C)^{-1}
d\zeta,
$$  
where $\Gamma$ is the boundary of a Cauchy domain $\Delta$ containing  
$\sigma(T_\C)$ in $U$. In fact, because $\sigma(T_\C)$ is conjugate symmetric, we may and shall assume 
that both $U$ and $\Gamma$ are conjugate symmetric. Because the
function $\zeta\mapsto F(\zeta)(\zeta-T_\C)^{-1}$ is analytic
in $U\setminus\sigma(T_\C)$, the integral does not depend on the
particular choice of the Cauchy domain $\Delta$ containing 
$\sigma(T_C)$.

A natural 
question is to find an appropriate condition to we have 
$F(T_\C)^\flat=F(T_\C)$, which would
imply the invariance of $\mathcal{V}$ under $F(T_\C)$. 
\end{Rem}

With the notation of Remark \ref{afcro}, we have the following.

% Thm 1

\begin{Thm}\label{afcro1} Let $U\subset\C$ be open and conjugate symmetric. If $F:U\mapsto\mathcal{B}(\mathcal{V}_\C)$ is analytic and $F(\zeta)^\flat=F(\bar{\zeta})$ for all $\zeta\in U$, then  $F(T_\C)^\flat=F(T_\C)$ for all $T\in\mathcal{B}(\mathcal{V})$ with $\sigma_\C(T)\subset U$. 
\end{Thm}

{\it Proof.}\, We use the  notation from  Remark \ref{afcro}, assuming, in
addition, that  $\Gamma$ is conjugate symmetric as well. We put 
$\Gamma_\pm:=\Gamma\cap\C_\pm$, where $\C_+$ (resp. $\C_-$) equals to $\{\lambda\in\C;\Im\lambda\ge0\}$ (resp. $\{\lambda\in\C;\Im\lambda\le0\}$). We write $\Gamma_+=\cup_{j=1}^m\Gamma_{j+}$, where $\Gamma_{j+}$ are the connected components
of $\Gamma_+$. Similarly, we write  $\Gamma_-=\cup_{j=1}^m\Gamma_{j-}$, where $\Gamma_{j-}$ are the connected components
of $\Gamma_-$, and $\Gamma_{j-}$ is the reflexion of 
$\Gamma_{j+}$ with respect of the real axis. 

As $\Gamma$ is a finite union of Jordan piecewise smooth closed curves, 
for each index
$j$ we have a parametrization $\phi_j:[0,1]\mapsto\C$, positively oriented, such that
$\phi_j([0,1])=\Gamma_{j+}$. Taking into account that the function $t\mapsto\overline{\phi_j(t)}$ is a 
parametrization of $\Gamma_{j-}$ negatively oriented, and setting 
$\Gamma_j=\Gamma_{j+}\cup\Gamma_{j-}$, we can write
$$
F_j(T_\C):=\frac{1}{2\pi i}\int_{\Gamma_j} F(\zeta)(\zeta-T_\C)^{-1}
d\zeta=
$$
$$
\frac{1}{2\pi i}\int_0^1 F(\phi_j(t))(\phi_j(t)-T_\C)^{-1}
\phi_j'(t)dt
$$
$$
-\frac{1}{2\pi i}\int_0^1 F(\overline{\phi_j(t)})(\overline{\phi_j(t)}-T_\C)^{-1}\overline{\phi_j'(t)}dt.
$$
Therefore,
$$
F_j(T_\C)^\flat=
-\frac{1}{2\pi i}\int_0^1 F(\phi_j(t))^\flat(\overline{\phi_j(t)}-T_\C)^{-1}\overline{\phi_j'(t)}dt
$$
$$
+\frac{1}{2\pi i}\int_0^1 F(\overline{\phi_j(t)})^\flat(\phi_j(t)-T_\C)^{-1}\phi_j'(t)dt.
$$
According to our assumption on the function $F$, we obtain
$F_j(T_\C)=F_j(T_\C)^\flat$ for all $j$, and therefore 
$$
F(T_\C)^\flat=\sum_{j=1}^mF_j(T_\C)^\flat=\sum_{j=1}^mF_j(T_\C)=F(T_\C), $$
which concludes the proof.

% Rem 3

\begin{Rem}\label{stem_anal}\rm If $\A$ is a unital real Banach algebra,  
$\A_\C$ its complexification, and $U\subset\C$ is open, we denote by 
$\mathcal{O}(U,\A_\C)$ the algebra of all analytic $\A_\C$-valued functions. If $U$ is conjugate symmetric, and 
$\A_\C\ni a\mapsto \bar{a}\in
\A_\C$ is its natural conjugation, we denote by $\mathcal{O}_s(U,\A_\C)$ the real subalgebra of $\mathcal{O}(U,\A_\C)$ consisting of those functions $F$ with the property $F(\bar{\zeta})=\overline{F(\zeta)}$
for all $\zeta\in U$. Adapting a well known terminology, such 
functions will be called ($\A_\C$-{\it valued $)$ stem functions}.

When $\A=\R$, so $\A_\C=\C$, the space $\mathcal{O}_s(U,\C)$ will be denoted by 
$\mathcal{O}_s(U)$, which is a real algebra. Note that 
$\mathcal{O}_s(U,\A_\C)$ is also a bilateral $\mathcal{O}_s(U)$-module. 
 
In the next result, we identify the algebra 
$\mathcal{B}(\mathcal{V})$ with a subalgebra of 
$\mathcal{B}(\mathcal{V})_\C$. In ths case, when $F\in \mathcal{O}_s(U,\mathcal{B}(\mathcal{V})_\C)$, we shall write
$$
F(T)=\frac{1}{2\pi i}\int_\Gamma F(\zeta)(\zeta-T)^{-1}
d\zeta,
$$  
noting that the right hand side of this formula belongs to 
$\mathcal{B}(\mathcal{V})$, by Theorem \ref{afcro1}.
\end{Rem} 

The properties of the map $F\mapsto F(T)$, which can be called
the {\it (left) analytic functional calculus of} $T$, are summarized by 
the following.

% Thm 2 

\begin{Thm}\label{R_afc} Let $\V$ be a real Banach space, let $U\subset\C$ be a conjugate symmetric open set, and let $T\in\mathcal{B}(\mathcal{V})$, with
$\sigma_\C(T)\subset U$.  Then the assignment 
$$
{\mathcal O}_s(U,\mathcal{B}(\mathcal{V})_\C)\ni F\mapsto F(T)\in\mathcal{B}(\mathcal{V})
$$ 
is an $\R$-linear map, and the map
$$
{\mathcal O}_s(U)\ni f\mapsto f(T)\in\mathcal{B}(\mathcal{V})
$$
is a unital real algebra morphism.

Moreover, the following properties are true:

(1) For all   $F\in\mathcal{O}_s(U,\mathcal{B}(\mathcal{V})_\C),\, f\in{\mathcal O}_s(U)$, we have $(Ff)(T)=F(T)f(T)$.

(2) For every polynomial $P(\zeta)=\sum_{n=0}^m A_n\zeta^n,\,\zeta\in\C$, with $A_n\in\mathcal{B}(\mathcal{V})$  for all $n=0,1,\ldots,m$, we have  $P(T)=\sum_{n=0}^m A_n T^n\in\mathcal{B}(\mathcal{V})$.  
\end{Thm}

{\it Proof.\,} The arguments are more or less standard (see
\cite{DuSc}). The $\R$-linearity of the maps
$$
{\mathcal O}_s(U,\mathcal{B}(\mathcal{V})_\C)\ni F\mapsto F(T)\in\mathcal{B}(\mathcal{V}),\,
{\mathcal O}_s(U)\ni f\mapsto f(T)\in\mathcal{B}(\mathcal{V}),
$$
is clear. The second one is actually multiplicative, which follows from the multiplicativiry of the usual analytic functional calculus  of $T$.

In fact, we have a more general property, specifically 
$(Ff)(T)=F(T)f(T)$ for all   $F\in\mathcal{O}_s(U,\mathcal{B}(\mathcal{V})_\C),\, f\in{\mathcal O}_s(U)$. This follows from the equalities,
$$
(Ff)(T)=\frac{1}{2\pi i}\int_{\Gamma_0} F(\zeta)f(\zeta)(\zeta-T)^{-1}d\zeta=
$$
$$
\left(\frac{1}{2\pi i}\int_{\Gamma_0} F(\zeta)(\zeta-T)^{-1}d\zeta\right)
\left(\frac{1}{2\pi i}\int_{\Gamma} f(\eta)(\eta-T)^{-1}d\eta\right)=F(T)f(T),
$$
obtained as in the classical case (see \cite{DuSc}, Section VII.3), which holds because $f$ is $\C$-valued and commutes with the operators in $\mathcal{B}(\mathcal{V})$. Here $\Gamma,\,\Gamma_0$ are the boundaries of two Cauchy
domains $\Delta,\,\Delta_0$ respectively, such that $\Delta\supset
\bar{\Delta}_0$, and $\Delta_0$ contains $\sigma(T)$. 

Note that, in particular,  for every polynomial $P(\zeta)=\sum_{n=0}^m A_n\zeta^n$ with $A_n\in\mathcal{B}(\mathcal{V})$  for all $n=0,1,\ldots,m$, we have  $P(T)=\sum_{n=0}^m A_n q^n\in\mathcal{B}(\mathcal{V})$ for all $T\in\mathcal{B}(\mathcal{V})$.

% Ex 1 
 
 \begin{Exa}\label{ex1}\rm Let $\mathcal{V}=\R^2$, so $\mathcal{V}_\C=\C^2$, endowed with its natural Hilbert space
structure. Let us first observe that we have 
$$
S=\left(\begin{array}{cc} a_1 & a_2 \\ a_3 &  a_4
\end{array}\right)\,\,\Longleftrightarrow
S^\flat=\left(\begin{array}{cc} \bar{a}_1 & \bar{a}_2 \\ 
\bar{a}_3 &  \bar{a}_4\end{array}\right),
$$
for all $a_1,a_2,a_3,a_4\in\C$. 

Next we consider the operator $T\in\mathcal{B}(\R^2)$ given by the matrix 
$$
T=\left(\begin{array}{cc} u & v \\ -v & u
\end{array}\right),
$$
where $u,v\in\R, v\neq0$. The extension $T_\C$ of the operator 
$T$ to $\C^2$, which is a normal operator, is given by the same formula. Note that
$$
\sigma_\C(T)=\{\lambda\in\C;(\lambda-u)^2+v^2=0\}=
\{u\pm iv\}=\sigma(T_\C).
$$
Note also that the vectors  $\nu_\pm=(\sqrt{2})^{-1}(1,\pm i)$ 
are normalized eigenvectors for $T_\C$ corresponding to the eigenvalues $u\pm iv$, respectively.
The spectral projections of $T_\C$ corresponding to 
these  eigenvalues are given by 
$$
E_\pm(T_\C){\bf w}=\langle {\bf w},\nu_\pm\rangle \nu_\pm=\frac{1}{2}\left(\begin{array}{cc} 1 & \mp i \\ \pm i & 1 \end{array}\right)\left(\begin{array}{c} w_1 \\ w_2 \end{array}\right),
$$ 
for all ${\bf w}=(w_1,w_2)\in\C^2$.

Let $U\subset\C$ be an open set with $U\supset\{u\pm iv\}$,
and let $F:U\mapsto\mathcal{B}(\C^2)$ be analytic. 
We shall compute explicitly $F(T_\C)$. Let $\Delta$ be a Cauchy 
domain contained in $U$  with its boundary $\Gamma$, and containing  the points $u\pm iv$. Assuming $v>0$, we have
$$
F(T_\C)=\frac{1}{2\pi i}\int_\Gamma F(\zeta)(\zeta-T_\C)^{-1}
d\zeta=
$$
$$
F(u+iv)E_+(T_\C)+F(u-iv)E_-(T_\C)=
$$
$$
\frac{1}{2}F(u+iv)\left(\begin{array}{cc} 1 & -i \\ i & 1 \end{array}\right)+\frac{1}{2}F(u-iv)\left(\begin{array}{cc} 1 & i \\ -i & 1 \end{array}\right).
$$
 
Assume now that $F(T_\C)^\flat=F(T_\C)$. Then we must have

$$
(F(u+iv)-F(u-iv)^\flat)\left(\begin{array}{cc} 1 & -i \\ i & 1
\end{array}\right)= 
(F(u+iv)^\flat-F(u-iv))\left(\begin{array}{cc} 1 & i \\ -i & 1
\end{array}\right).
$$
We also have the equalities
$$
\left(\begin{array}{cc} 1 & -i \\ i & 1
\end{array}\right)\left(\begin{array}{c} 1 \\ i \end{array} \right)=2\left(\begin{array}{c} 1 \\ i \end{array}\right),\,\,
\left(\begin{array}{cc} 1 & -i \\ i & 1
\end{array}\right)\left(\begin{array}{c} 1 \\ -i \end{array} \right)=0,
$$
$$
\left(\begin{array}{cc} 1 & i \\ -i & 1
\end{array}\right)\left(\begin{array}{c} 1 \\ -i \end{array} \right)=2\left(\begin{array}{c} 1 \\ -i \end{array}\right),\,\,
\left(\begin{array}{cc} 1 & i \\ -i & 1
\end{array}\right)\left(\begin{array}{c} 1 \\ i \end{array} \right)=0,
$$
Using these equalities, we finally deduce that
$$
(F(u+iv)-F(u-iv)^\flat)\left(\begin{array}{c} 1 \\ i \end{array}\right)=0,
$$
and
$$
(F(u-iv)-F(u+iv)^\flat)\left(\begin{array}{c} 1 \\ -i \end{array}\right)=0,
$$
which are necessary conditions for the equality $F(T_\C)^\flat=F(T_\C)$. 
As a matter of fact, this example shows, in particular, that the condition
$F(\zeta)^\flat=F(\bar{\zeta})$ for all $\zeta\in U$,
used in Theorem \ref{afcro1}, is sufficient but it might not be always necessary. 
\end{Exa}

\section{Analytic Functional Calculus for Quaternionic 
\\ Operators}
 
\subsection{Quaternionic Spectrum}

We now recall some  known definitions and elementary
facts (see, for instance, \cite{CoSaSt}, Section 4.6, and/or
\cite{Vas5}). 

Let $\H$ be the abstract algebra of quaternions, which is the four-dimensional $\R$-algebra with 
unit $1$, generated by the ''imaginary units`` $\{\bf{j,k,l}\}$,  
which satisfy
$$
{\bf jk=-kj=l,\,kl=-lk=j,\,lj=-jl=k,\,jj=kk=ll}=-1.
$$

We may assume that $\H\supset\R$ identifying every  number
$x\in\R$ with the element $x1\in\H$.

The algebra $\H$ has a natural multiplicative norm given by
$$
\Vert {\bf x}\Vert=\sqrt{x_0^2+x_1^2+x_2^2+x_0^2},\,\,{\bf x}= x_0+x_1{\bf j}+x_2{\bf k}+x_3{\bf l},\,\,x_0,x_1,x_2,x_3\in\R,
$$
and a natural involution
$$
\H\ni{\bf x} = x_0+x_1{\bf j}+x_2{\bf k}+x_3{\bf l}\mapsto
{\bf x}^*= x_0-x_1{\bf j}-x_2{\bf k}-x_3{\bf l}\in\H.
$$
Note that ${\bf x}{\bf x}^*={\bf x}^*{\bf x}=\Vert{\bf x}\Vert^2$, implying, in particular, that every element ${\bf x}\in\H\setminus\{0\}$ is invertible, and ${\bf x}^{-1}=
\Vert {\bf x}\Vert^{-2}{\bf x}^*$.

For an arbitrary quaternion ${\bf x}= x_0+x_1{\bf j}+x_2{\bf k}+x_3{\bf l},\,\,x_0,x_1,x_2,x_3\in\R$, we set $\Re{\bf x}=x_0=
({\bf x}+{\bf x}^*)/2$, and $\Im{\bf x}=x_1{\bf j}+x_2{\bf k}+x_3{\bf l}=({\bf x}-{\bf x}^*)/2$, that is, the {\it real} and 
{\it imaginary part} of ${\bf x}$, respectively.
 
We consider the complexification $\C\otimes_\R\H$
of the $\R$-algebra $\H$ (see also \cite{GhMoPe}), which will be
identified with the  direct sum $\M=\H+i\H$. 
Of course, the algebra $\M$ contains the complex field $\C$.  Moreover, in the algebra $\M$, the elements of $\H$ commute with all complex numbers.  In particular, the ''imaginary units`` 
$\bf j,k,l$ of the algebra $\H$  are  independent of and commute with the  imaginary unit $i$ of the complex plane $\C$.

In the algebra $\M$, there also exists a natural conjugation given
by $\bar{\bf a}={\bf b}-i{\bf c}$, where ${\bf a}={\bf b}+i{\bf c}$ is arbitrary in $\M$, with ${\bf b},{\bf c}\in\H$ (see also
\cite{GhMoPe}). Note that $\overline{\bf a+b}=\bar{\bf a}+\bar{\bf b}$, and  $\overline{\bf ab}=\bar{\bf a}\bar{\bf b}$, in particular  $\overline{r\bf a}=r\bar{\bf a}$ for all ${\bf a},{\bf b}\in\M$, and $r\in\R$.  Moreover, $\bar{{\bf a}}={\bf a}$ if and only if ${\bf a}\in\H$, which is a useful characterization of the elements of $\H$ among those of $\M$. 
 
% Rem 4

\begin{Rem}\label{Qspectrum}\rm
In the algebra $\M$ we have the identities
$$
(\lambda-{\bf x}^*)(\lambda-{\bf x})=(\lambda-{\bf x})(\lambda-{\bf x}^*)=\lambda^2-
\lambda({\bf x}+{\bf x}^*)+\Vert {\bf x}\Vert^2\in\C,
$$
for all $\lambda\in\C$ and ${\bf x}\in\H$. If the complex number 
$\lambda^2-2\lambda\Re{\bf x}+\Vert {\bf x}\Vert^2$ is 
nonnull, then both element $\lambda-{\bf x}^*,\,\lambda-{\bf x}$ are invertible. Conversely, if $\lambda-{\bf x}$ is invertible, we must have 
$\lambda^2-2\lambda\Re{\bf x}+\Vert {\bf x}\Vert^2$ nonnull; otherwise we 
would have $\lambda={\bf x}^*\in\R$, so $\lambda={\bf x}\in\R$, which is not possible. Therefore, the element
$\lambda-{\bf x}\in\M$ is invertible if and only if the complex number $\lambda^2-2\lambda\Re{\bf x}+\Vert {\bf x}\Vert^2$ is 
nonnull. Hence, the element      
$\lambda-{\bf x}\in\M$ is not  invertible if and only if $\lambda=
\Re{\bf x}\pm i\Vert\Im{\bf x}\Vert$. In this way, the {\it spectrum} of a quaternion ${\bf x}\in\H$ is given by the equality  $\sigma({\bf x})=\{s_\pm(\bf x)\}$, where 
$s_\pm(\bf x)=\Re{\bf x}\pm i\Vert\Im{\bf x}\Vert$ are the {\it eigenvalues} of $\bf x$ (see also \cite{Vas4,Vas5}).  

The polynomial $P_{\bf x}(\lambda)=\lambda^2-2\lambda\Re{\bf x}+\Vert {\bf x}\Vert^2$ is the {\it minimal polynomial} of $\bf x$.  In fact,
the equality $\sigma({\bf y})= \sigma({\bf x})$ for some ${\bf x,y}\in\H$
is an equivalence relation in the algebra $\H$, which holds if and only if  $P_{\bf x}=P_{\bf y}$. In fact, setting $\mathbb{S}=\{\mathfrak{\kappa}\in\H;\Re\mathfrak{\kappa}=0, \Vert \mathfrak{\kappa}\Vert=1\}$ (that is the unit sphere of purely imaginary quaternions), representig an arbitrary
quaternion $\bf x$ under the form $x_0+y_0 \mathfrak{\kappa}_0$, with 
$x_0,y_0\in\R$ and $\mathfrak{\kappa}_0\in\mathbb{S}$, a quaternion $\bf y$
is equivalent to $\bf x$ if anf only if it is of the form $x_0+y_0 \mathfrak{\kappa}$ for some $\mathfrak{\kappa}\in\mathbb{S}$ (see \cite{Bre} or \cite{Vas5} for some details). 
\end{Rem}

% Rem 5

\begin{Rem}\label{Hspace}\rm
Following \cite{CoSaSt}, a {\it right $\H$-vector space} 
$\mathcal{V}$ is a real vector space 
having a right multiplication with the elements of $\H$, such that $(x+y){\bf q}=x{\bf q}+y{\bf q},\,x({\bf q}+{\bf s})=
x{\bf q}+x{\bf s},\, x({\bf q}{\bf s})=(x{\bf q}){\bf s}$
for all $x,y\in\mathcal{V}$ and ${\bf q},{\bf s}\in\H$.

 If $\mathcal{V}$ is also a Banach space the operator 
$T\in\mathcal{B(V)}$ is {\it right $\H$-linear} if 
$T(x{\bf q})=T(x){\bf q}$ for all $x\in\mathcal{V}$ and
$\q\in\H$. The set of right $\H$ linear operators will be 
denoted by $\mathcal{B^{\rm r}(V)}$, which is, in particular, a unital 
real algebra.

In a similar way, one defines the concept of a {\it left $\H$-vector space}. A real vector space $\mathcal{V}$ will be said to be an {\it $\H$-vector space} if it is simultaneously a right $\H$- and a left $\H$-vector space. As noticed in \cite{CoSaSt},
it is the framework of  $\H$-vector spaces an appropriate one
for the study of right $\H$-linear operators.

If $\V$ is $\H$-vector space which is also a Banach space, then
$\V$ is said to be a {\it Banach $\H$-space}. In this case, we also assume that $ R_\q\in \mathcal{B}(\V)$, and 
the map $\H\ni\q\mapsto R_\q\in 
\mathcal{B}(\V)$ is norm continuous, where $R_{\bf q}$ is the right multiplication of the elements of 
$\mathcal{V}$ by a given quaternion ${\bf q}\in\H$. Similarly, 
if  $L_{\bf q}$ is the left multiplication of the elements of 
$\mathcal{V}$ by the quaternion ${\bf q}\in\H$, we assume that
 $ L_\q\in \mathcal{B}(\V)$ for all ${\bf q}\in\H$, and that 
the map $\H\ni\q\mapsto L_\q\in \mathcal{B}(\V)$ is norm continuous. Note also that 
$$
\mathcal{B^{\rm r}(V)}=\{T\in\mathcal{B(V)};TR_\q=R_\q T,\,\q\in\H\}.
$$

To adapt the discussion regarding the real algebras to this case, we first consider the complexification $\V_\C$ of $\V$. Because $\V$ is an $\H$-bimodule, the space $\V_\C$ is actually an $\M$-bimodule, via the multiplications
$$
(\q+i\s)(x+iy)=\q x-\s y+i(\q y+\s x),
(x+iy)(\q+i\s)=x\q-y\s+i(y\q+x\s),
$$ 
for all $\q+i\s\in\M,\,\q,\s\in\H,\,x+iy\in\V_\C,\,x,y\in\V$. Moreover, the operator 
$T_\C$ is right $\M$-linear, that is $T_\C((x+iy)(\q+i\s))=
T_\C(x+iy)(\q+i\s)$ for all $\q+i\s\in\M,\,x+iy\in\V_\C$, via a direct computation.

Let $C$ be the conjugation of $\V_\C$. As in the real case, for every
$S\in \mathcal{B}(\V_\C)$, we put $S^\flat=CSC$. The left and  right
multiplication with the quaternion $\q$ on $\V_\C$ will be also denoted by $L_\q,R_\q$, respectively, as elements of $\mathcal{B}(\V_\C)$. We set
$$
\mathcal{B}^{\rm r}(\V_\C)=\{S\in \mathcal{B}(\V_\C); SR_\q=R_\q S,\,\q\in\H\},
$$
which is a unital complex algebra containing all operators $L_\q,
\q\in\H$.
Note that if $S\in\mathcal{B}^{\rm r}(\V_\C)$, then $S^\flat\in\mathcal{B}^{\rm r}(\V_\C)$. Indeed, because $CR_\q=R_\q C$, we also have  $S^\flat R_\q=R_\q S^\flat$. In fact, as we have
$(S+S^\flat)(\V)\subset\V$ and $i(S-S^\flat)(\V)\subset\V$, it folows that  
the algebras $\mathcal{B}^{\rm r}(\V_\C),\,\mathcal{B}^{\rm r}(\V)_\C$
are isomorphic, and they will be often identified, where 
$\mathcal{B^{\rm r}(V)}_\C=\mathcal{B^{\rm r}(V)}+i\mathcal{B^{\rm r}(V)}$ is the complexification of $\mathcal{B^{\rm r}(V)}$, which is also a unital
complex Banach algebra.
\end{Rem}

Looking at the Definition 4.8.1 from \cite{CoSaSt}
(see also \cite{CoGaKi}), we give the folowing. 

% Def 1

\begin{Def}\label{Q-spectrum}\rm For a given operator $T\in\mathcal{B^{\rm r}(V)}$, the set 
$$
\sigma_\H(T):=\{\q\in\H; T^2-2(\Re\q)T+\Vert\q\Vert^2)\,\,
{\rm not}\,\,{\rm invertible}\}
$$
is called the {\it quaternionic spectrum} (or simply the $Q$-{\it spectrum}) of $T$.
 
The complement $\rho_\H(T)=\H\setminus\sigma_\H(T)$ is called 
the {\it quaternionic resolvent} (or simply the $Q$-{\it resolvent}) of $T$.
\end{Def}

 Note that, if $\q\in\sigma_\H(T)$), then 
 $\{{\bf s}\in\H;\sigma({\bf s})=\sigma(\q)\}\subset\sigma_\H(T)$.
  
 Assuming that $\V$ is a  Banach $\H$-space, then
$\mathcal{B^{\rm r}(V)}$ is a unital real Banach $\H$-algebra 
(that is, a Banach algebra which also a Banach $\H$-space), via the 
algebraic operations  $(\q T)(x)=\q T(x)$, and $(T\q)(x)=T(\q x)$ for all
$\q\in\H$ and $x\in\V$. Hence the complexification $\mathcal{B^{\rm r}(V)}_\C$ is, in particular, a unital complex Banach algebra. Also note that the  complex numbers,
regarded as elements of $\mathcal{B^{\rm r}(V)}_\C$, commute with
the elements of $\mathcal{B^{\rm r}(V)}$. For this reason, for
each $T\in\mathcal{B^{\rm r}(V)}$ we have the resolvent set
$$
\rho_\C(T)=\{\lambda\in\C;(T^2-2(\Re\lambda)T+\vert\lambda\vert^2)^{-1}\in\mathcal{B^{\rm r}(V)}\}= 
$$
$$
\{\lambda\in\C;(\lambda-T_\C)^{-1}\in \mathcal{B^{\rm r}(V}_\C)\}=\rho(T_\C),
$$  
and the associated spectrum $\sigma_\C(T)=\sigma(T_\C)$. 
 
Clearly, there exists a strong connexion between $\sigma_\H(T)$ and 
$\sigma_\C(T)$. In fact, the set $\sigma_\C(T)$ looks like a ''complex border`` of the set $\sigma_\H(T)$. Specifically, we can prove the following.

% Lem 2

\begin{Lem}\label{spec_eg0} For every $T\in\mathcal{B^{\rm r}(V)}$ we have the
equalities 
\begin{equation}\label{spec_eg} 
\sigma_\H(T)=\{\q\in\H;\sigma_\C(T)\cap\sigma(\q)\neq\emptyset\}.
\end{equation}
and 
\begin{equation}\label{spec_eg1}
\sigma_\C(T)=\{\lambda\in\sigma(\q);\q\in\sigma_\H(T)\}.
\end{equation}
 
\end{Lem}

{\it Proof.} Let us prove (\ref{spec_eg}). If $\q\in \sigma_\H(T)$, and so the $T^2-2(\Re\q)T+\Vert\q\Vert^2$ is not invertible, choosing  $\lambda\in\{\Re\q\pm i\Vert\Im\q\Vert\}=\sigma(\q)$, we clearly have $T^2-2(\Re\lambda)T+\vert\lambda\vert^2$ not invertible, implying $\lambda\in\sigma_\C(T)\cap\sigma(\q)\neq\emptyset$. 

Conversely, if for some $\q\in\H$ there exists $\lambda\in\sigma_\C(T)\cap\sigma(\q)$, and so $T^2-2(\Re\lambda)T+\vert\lambda\vert^2=T^2-2(\Re\q)T+\Vert\q\Vert^2$ is not invertible, implying $\q\in\sigma_\H(T)$. 

 We now prove (\ref{spec_eg1}). Let $\lambda\in\sigma_\C(T)$, so the operator $T^2-2(\Re\lambda)T+\vert\lambda\vert^2$ is not invertible. Setting
$\q=\Re(\lambda)+\Vert\Im\lambda\Vert\kappa$, with $\kappa\in\mathbb{S}$, we have 
$\lambda\in\sigma(\q)$. Moreover, $T^2+2(\Re\q)T+\Vert\q\Vert^2$ is not 
invertible, and so $\q\in\sigma_ \H(T)$. 

Conversely, if $\lambda\in\sigma(\q)$ for some $\q\in\sigma_\H(T)$, then
$\lambda\in\{\Re\q\pm i\Vert\Im(\q)\Vert\}$, showing that
$T^2-2\Re(\lambda)T+\vert\lambda\vert^2=T^2+2(\Re\q)T+\Vert\q\Vert^2$ is
not invertible. 
\medskip

{\bf Remark} As expected, the set $\sigma_\H(T)$ is nonempty and bounded, which follows easily from Lemma \ref{spec_eg0}. It is also compact, as a consequence of Definition \ref{Q-spectrum}, because the set of invertible elements in $\mathcal{B^{\rm r}(V)}$ is open.

We recall that a subset $\Omega\subset\H$ is said to be {\it spectrally saturated} 
(see \cite{Vas4},\cite{Vas5}) if whenever $\sigma({\bf h})=\sigma({\bf q})$ for some ${\bf h}\in\H$ and ${\bf q}\in\Omega$, we also have ${\bf h}\in\Omega$. As noticed in \cite{Vas4} and \cite{Vas5}, this concept coincides with that of {\it axially symmetric set}, introduced in \cite{CoSaSt}.  

Note that the subset $\sigma_\H(T)$ spectrally saturated. 

\subsection{Analytic Functional Calculus}

If $\V$ is a Banach $\H$-space, because 
$\mathcal{B^{\rm r}(\V)}$
is real Banach space, each operator $T\in\mathcal{B^{\rm r}(\V)}$ has a complex spectrum $\sigma_\C(T)$. Therefore, applying the
corresponding result for real operators, we may construct
an analytic functional calculus using the classical Riesz-Dunford functional calculus, in a slightly  generalized form. In this case, our basic complex algebra is $\mathcal{B}^{\rm r}(\mathcal{V})_\C$, endowed with the  conjugation
$\mathcal{B}^{\rm r}(\mathcal{V})_\C\ni S\mapsto S^\flat\in\mathcal{B}^{\rm r}(\mathcal{V})_\C$.
 
% Thm 3

\begin{Thm}\label{afcro2}  Let $U\subset\C$ be open and conjugate symmetric. If $F:U\mapsto\mathcal{B}^{\rm r}(\mathcal{V}_\C)$ is analytic and $F(\zeta)^\flat=F(\bar{\zeta})$ for all $\zeta\in U$, then  $F(T_\C)^\flat=F(T_\C)$ for all $T\in\mathcal{B}^{\rm r}(\mathcal{V})$ with 
$\sigma_\C(T)\subset U$.
\end{Thm}

Both the statement and the proof of Theorem \ref{afcro2} are similar to 
those of Theorem \ref{afcro1}, and will be omitted.

As in the real case, we may identify the algebra 
$\mathcal{B}^{\rm r}(\mathcal{V})$ with a subalgebra of 
$\mathcal{B}^{\rm r}(\mathcal{V})_\C$. In ths case, when $F\in \mathcal{O}_s(U,\mathcal{B}^{\rm r}(\mathcal{V})_\C)=
\{F\in{\mathcal O}(U,\mathcal{B}^{\rm r}(\V)_\C);F(\bar{\zeta})=F(\zeta)^{\flat}\,\, \forall \zeta\in U\}$ (see also Remark \ref{spec_lm}), we can write, via the previous Theorem,
$$
F(T)=\frac{1}{2\pi i}\int_\Gamma F(\zeta)(\zeta-T)^{-1}
d\zeta\in\mathcal{B}^{\rm r}(\mathcal{V}),
$$  
for a suitable choice of $\Gamma$. 

The next result provides an {\it analytic functional calculus} for 
operators from the real algebra $\mathcal{B}^{\rm r}(\mathcal{V})$. 

% Thm 4 
\begin{Thm}\label{H_afc} Let $\V$ be a Banach $\H$-space, let $U\subset\C$ be a conjugate symmetric open set, and let $T\in\mathcal{B}^{\rm r}(\mathcal{V})$, with
$\sigma_\C(T)\subset U$.  Then the map 
$$
{\mathcal O}_s(U,\mathcal{B}^{\rm r}(\mathcal{V})_\C)\ni F\mapsto F(T)\in\mathcal{B}^{\rm r}(\mathcal{V})
$$ 
is $\R$-linear, and the map
$$
{\mathcal O}_s(U)\ni f\mapsto f(T)\in\mathcal{B}^{\rm r}(\mathcal{V})
$$
is a unital real algebra morphism.

Moreover, the following properties are true:

$(1)$ For all   $F\in\mathcal{O}_s(U,\mathcal{B}^{\rm r}(\mathcal{V})_\C),\, f\in{\mathcal O}_s(U)$, we have $(Ff)(T)=F(T)f(T)$.

$(2)$ For every polynomial $P(\zeta)=\sum_{n=0}^m A_n\zeta^n,\,\zeta\in\C$, with $A_n\in\mathcal{B}^{\rm r}(\mathcal{V})$  for all $n=0,1,\ldots,m$, we have  $P(T)=\sum_{n=0}^m A_n T^n\in\mathcal{B}^{\rm r}(\mathcal{V})$.  
\end{Thm} 

The proof of this result is similar to that of Theorem \ref{R_afc} and will be omitted.

% Rem 6

\begin{Rem}\label{leftmult}\rm The algebra $\H$ is, in particular, a Banach
$\H$-space. As already noticed, the left multiplications $L_\q,\,\q\in\H,$ are elements of 
$\mathcal{B}^{\rm r}(\H)$. In fact, the map $\H\ni\q\mapsto L_\q\in \mathcal{B}^{\rm r}(\H)$ is a injective morphism of real
algebras allowing the identification of $\H$ with a subalgebra
of $\mathcal{B}^{\rm r}(\H)$.  
\end{Rem}

Let $\Omega\subset\H$ be a spectrally saturated open set, and 
let $U=\mathfrak{S}(\Omega):=\{\lambda\in\C, \exists \q\in\Omega, \lambda\in\sigma(\q)\}$, which is open and conjugate symmetric (see \cite{Vas5}).  Denotig by $f_\H$ the function $\Omega\ni\q\mapsto
f(\q),\q\in\Omega$, for every $f\in\mathcal{O}_s(U)$, we set  
$$
\mathcal{R}(\Omega):=\{f_\H; f\in\mathcal{O}_s(U)\},
$$
which is a commutative real algebra. Defining the function $F_\H$ in a similar way for each $F\in\mathcal{O}_s(U,\M)$, we set
$$
\mathcal{R}(\Omega,\H):=\{F_\H; F\in\mathcal{O}_s(U,\M)\},
$$
which, according to the next theorem, is a right $\mathcal{R}(\Omega)$-module. 

The next result is an analytic functional calculus for quaternions (see \cite{Vas5}, Theorem 5), obtained as a particular case of Theorem  \ref{H_afc} (see also its predecessor in \cite{CoSaSt}).

% Thm 5

\begin{Thm}\label{H_afc0} Let $\Omega\subset\H$ be a spectrally
saturated open set, and let $U=\mathfrak{S}(\Omega)$.
The space $\mathcal{R}(\Omega)$ is a unital commutative $\R$-algebra, the space $\mathcal{R}(\Omega,\H)$ is a right 
$\mathcal{R}(\Omega)$-module, the map
$$
{\mathcal O}_s(U,\M)\ni F\mapsto F_\H\in\mathcal{R}(\Omega,\H)
$$ 
is a right module isomorphism, and  its restriction
$$
{\mathcal O}_s(U)\ni f\mapsto f_\H\in\mathcal{R}(\Omega)
$$ 
is an $\R$-algebra isomorphism.

Moreover, for every polynomial\,\,$P(\zeta)=\sum_{n=0}^m a_n\zeta^n,\,\zeta\in\C$, with $a_n\in\H$  for all $n=0,1,\ldots,m$, we have  $P_\H(q)=\sum_{n=0}^m a_n q^n\in\H$ for all $q\in\H$.   
\end{Thm}

Most of the assertions of Theorem \ref{H_afc0} can be obtained directly 
from Theorem \ref{H_afc}. The injectivity of the map 
${\mathcal O}_s(U)\ni f\mapsto f_\H\in\mathcal{R}(\Omega)$, as well as an 
alternative complete proof, can be obtained as in the proof of Theorem 5 from \cite{Vas5}.

% Rem 7

\begin{Rem}\label{gen_afc}\rm That Theorems \ref{afcro2} and \ref{H_afc}
have practically the same proof as Theorems \ref{afcro1} and \ref{R_afc}
(respectively) is due to the fact that all of them can be obtained as particular
cases of more general results. Indeed, considering a unital real Banach 
algebra $\A$,  and its complexification $\A_\C$, identifying $\A$ with 
a real subalgebra of $\A_\C$, for a function $F\in\mathcal{O}_s(U,A_\C)$,
where $U\subset\C$ is open and conjugate symmetric,
the element $F(b)\in\A$ for each $b\in\A$ with $\sigma_\C(b)\subset U$. 
The assertion follows as in the proof of Theorem \ref{afcro1}. The other 
results also have their counterparts. We omit the details.  
\end{Rem}

% Rem 8

\begin{Rem}\label{twofc}\rm The space $\mathcal{R}(\Omega,\H)$ can be independently defined, and it consists of the set of all 
$\H$-valued functions, which are {\it slice regular} in the 
sense of \cite{CoSaSt}, Definition 4.1.1. They are used in  
\cite{CoSaSt} to define  a quaternionic  functional calculus
for quaternionic linear operators (see also
 \cite{CoGaKi}). Roughly speaking, given a quaternionic linear
 operator, each regular quaternionic-valued function  defined in a neighborhood $\Omega$ of its quaternionic spectrum is associated with
another quaternionic linear operator, replacing formally 
the quaternionic variable with that operator. This constraction 
is largely explained in the fourth chapter of  \cite{CoSaSt}.

Our Theorem \ref{H_afc} constructs an analytic functional calculus with functions from ${\mathcal O}_s(U,\mathcal{B}^{\rm r}(\mathcal{V})_\C)$, where $U$ is a a neighborhood of the complex spectrum  of a given quaternionic linear operator, leading to another quaternionic linear operator, replacing formally  the complex variable with that operator. We can show that those functional calculi are equivalent. This is a consequence of the 
fact that the class of regular quaternionic-valued function used by the construction in  \cite{CoSaSt} is isomorphic to  the class of analytic functions used in our Theorem \ref{H_afc0}. The advantage of our approach is its simplicity and a stronger connection with the classical approach, using spectra defined in the complex plane, and Cauchy type kernels partially commutative. 

Let us give a direct argument concerning the equivalence of those analytic functional calculi.  For an operator $T\in\mathcal{B}^{\rm r}(\mathcal{V})$, the so-called {\it right $S$-resolvent} is defined via the formula

\begin{equation}\label{kqfc} 
S_R^{-1}(\s,T)=-(T-\s^*)(T^2-2\Re(\s)T+\Vert\s\Vert)^{-1},\,\,\s\in
\rho_\H(T)
\end{equation}
(see \cite{CoSaSt}, formula (4.27)). Fixing an element $\kappa\in\mathbb{S}$, and a spectrally saturated open set $\Omega\subset\H$, for $\Phi\in\mathcal{R}(\Omega,\H)$ one sets 

\begin{equation}\label{qfc}
\Phi(T)=\frac{1}{2\pi}\int_{\partial(\Sigma_\kappa)}\Phi(\s)d\s_\kappa S_R^{-1}(\s,T),
\end{equation}
where $\Sigma\subset\Omega$ is a spectrally saturated open set containing 
$\sigma_\H(T)$, such that $\Sigma_\kappa=\{u+v\kappa\in\Sigma;u,v\in\R\}$
is a subset  whose boundary $\partial(\Sigma_\kappa)$ consists of a finite family of closed curves, piecewise smooth, positively oriented, and $d\s_\kappa=-\kappa du\wedge dv$. Formula (\ref{qfc}) is a (right) quaternionic functional calculus, as defined in \cite{CoSaSt}, Section 4.10.

Because the space $\mathcal{V}_\C$ is also an $\H$-space, we may extend these formulas to the 
operator $T_\C\in\mathcal{B}^{\rm r}(\mathcal{V}_\C)$, extending the operator $T$ to $T_\C$, and replacing $T$ by $T_\C$ in formulas (\ref{kqfc}) and (\ref{qfc}). For the function 
 $\Phi\in\mathcal{R}(\Omega,\H)$ there exists a function $F\in{\mathcal O}_s(U,\mathcal{B}^{\rm r}(\mathcal{V}_\C))$ such that $F_\H=\Phi$. Denoting by
 $\Gamma_\kappa$ the boundary of a Cauchy domain in $\C$ containing the compact set
 $\cup\{\sigma(\s);\s\in\overline{\Sigma_\kappa}\}$, we can write 
 $$
 \Phi(T_\C)=\frac{1}{2\pi}\int_{\partial(\Sigma_\kappa)}\left(\frac{1}{2\pi i}\int_{\Gamma_\kappa}F(\zeta)(\zeta-\s)^{-1}d\zeta\right)
 d\s_\kappa S_R^{-1}(\s,T_\C)=
 $$
 $$
 \frac{1}{2\pi i}\int_{\Gamma_\kappa}F(\zeta)\left(\frac{1}{2\pi}\int_{\partial(\Sigma_\kappa)}(\zeta-\s)^{-1}d\s_\kappa S_R^{-1}(\s,T_\C)\right)d\zeta.
 $$

It follows from the complex linearity of $S_R^{-1}(\s,T_\C)$, and from 
formula (4.49) in \cite{CoSaSt}, that
$$
(\zeta-\s)S_R^{-1}(\s,T_\C)=S_R^{-1}(\s,T_\C)(\zeta-T_\C)-1,
$$
whence
$$
(\zeta-\s)^{-1}S_R^{-1}(\s,T_\C)=S_R^{-1}(\s,T_\C)(\zeta-T_\C)^{-1}+
(\zeta-\s)^{-1}(\zeta-T_\C)^{-1},
$$      
and therefore,
 
$$
\frac{1}{2\pi}\int_{\partial(\Sigma_\kappa)}(\zeta-\s)^{-1}d\s_\kappa S_R^{-1}(\s,T_\C)=
 \frac{1}{2\pi}\int_{\partial(\Sigma_\kappa)}d\s_\kappa S_R^{-1}(\s,T_\C)
(\zeta-T_\C)^{-1}+
$$
$$
\frac{1}{2\pi}\int_{\partial(\Sigma_\kappa)}(\zeta-\s)^{-1}d\s_\kappa  
(\zeta-T_\C)^{-1}= (\zeta-T_\C)^{-1},
$$
because
$$
 \frac{1}{2\pi}\int_{\partial(\Sigma_\kappa)}d\s_\kappa S_R^{-1}(\s,T_\C)=1\,\,\,{\rm and}
 \,\,\,\frac{1}{2\pi}\int_{\partial(\Sigma_\kappa)}(\zeta-\s)^{-1}d\s_\kappa =0,
$$
as in Theorem 4.8.11 from \cite{CoSaSt}, since the $\M$-valued function $\s\mapsto(\zeta-\s)^{-1}$ is analytic in a neighborhood of the 
set $\overline{\Sigma_\kappa}\subset \C_\kappa$ for each $\zeta\in\Gamma_\kappa$, respectively. Therefore $\Phi(T_\C)=\Phi(T)_\C=F(T_\C)=F(T)_\C$, implying $\Phi(T)=F(T)$.

\end{Rem}

\section{Some Examples}

  % Ex 2
\begin{Exa}\label{ex2}\rm One of the simplest Banach $\H$-space
is the space $\H$ itself. As already noticed (see Remark \ref{leftmult}), taking $\V=\H$, so  $\V_\C=\M$, and fixing an element $\q\in\H$, we may consider the operator $L_\q\in\mathcal{B}^{\rm r}(\H)$, whose complex spectrum is given by
$\sigma_\C(L_\q)=\sigma(\q)=\{\Re\q\pm i\Vert\Im\q\Vert\}$.  If 
$U\subset\C$ is conjugate symmetric open set containing 
$\sigma_\C(L_\q)$, and $F\in\mathcal{O}_s(U,\M)$, then we have

\begin{equation}\label{GFC2}
F({L_\q})=F(s_+({\bf q}))\iota_+(\mathfrak{s}_{\tilde{\bf q}})+F(s_-({\bf q}))\iota_-(\mathfrak{s}_{\tilde{\bf q}}) \in\M,
\end{equation} 
where $s_\pm(\q)=\Re\q\pm i\Vert\Im\q\Vert$, 
$\tilde{\bf q}=\Im\bf q,\,\mathfrak{s}_{\tilde{\bf q}}=\tilde{\bf q}\Vert\tilde{\bf q}\Vert^{-1 }$, and  
$\iota_\pm(\mathfrak{s}_{\tilde{\bf q}})=2^{-1}
(1\mp i\mathfrak{s}_{\tilde{\bf q}})$ (see \cite{Vas5}, Remark 3).
\end{Exa}
 
% Ex 3 

\begin{Exa}\label{spec_lm}\rm Let $\X$ be a topological compact space, and let $C(\X,\M)$ be the space of $\M$-valued continuous
functions on $\X$. Then $C(\X,\H)$ is the real subspace of  
$C(\X,\M)$ consisting of $\H$-valued functions, which is also 
a Banach $\H$-space with respect to the operations $(\q F)(x)=
\q F(x)$ and $(F\q)(x)=F(x)\q$ for all $F\in C(\X,\H)$ and $x\in\X$. Moreover,  $C(\X,\H)_\C=C(\X,\H_\C)=C(\X,\M)$. 

We fix a function $\Theta\in C(\X,\H)$ and define the operator
$T\in\mathcal{B}(C(\X,\H))$ by the relation $(TF)(x)=\Theta(x)F(x)$ for all $F\in C(\X,\H)$ and $x\in\X$. Note that
$(T(F\q))(x)=\Theta(x)F(x)\q=((TF)\q)(x)$ for all 
$F\in C(\X,\H),\q\in\H$,  and $x\in\X$. In othe words, 
$T\in\mathcal{B}^{\rm r}(C(\X,\H))$. Note also that the operator $T$
is invertible if and only if the function $\Theta$ has no zero in $\X$. 

Let us compute the $Q$-spectrum of $T$. According to Definition \ref{Q-spectrum}, we have
$$
\rho_\H(T)=\{\q\in\H; (T^2-2\Re \q\,T+\Vert \q\Vert^2)^{-1}\in\mathcal{B}^{\rm r}(C(\X,\H))\}.
$$

Consequently, $\q\in\sigma_\H(T)$ if and only if zero is in the range of the function 
$$
\tau(\q,x):=\Theta(x)^2-2\Re \q\,\Theta(x)+\Vert \q\Vert^2
,\, x\in\mathfrak{X}.
$$ 
Similarly,

$$
\rho_\C(T)=\{\lambda\in\C; (T^2-2\Re \lambda\,T+\Vert \lambda\Vert^2)^{-1}\in\mathcal{B}^{\rm r}(C(\X,\H))\},
$$
and so $\lambda\in \sigma_\C(T)$ if and only if zero is in the range of the function 
$$
\tau(\lambda,x):=\Theta(x)^2-2\Re \lambda\,\Theta(x)+\vert 
\lambda\vert^2,\,x\in\mathfrak{X}.
$$ 

Looking for solutions $u+iv,u,v\in\R$, of the equation
$(u-\Theta(x))^2+v^2=0$,  a direct calculation shows that
$u=\Re\Theta(x)$ and $v=\pm\Vert\Im\Theta(x)\Vert$. Hence
$$
\sigma_\C(T)=\{\Re\Theta(x)\pm i\Vert\Im\Theta(x)\Vert;x\in\mathfrak{X}\}=\cup_{x\in\mathfrak{X}}\sigma(\Theta(x)).
$$
 
 Of course,  for every open conjugate symmetric subset
$U\subset\C$ containing $\sigma_\C(T)$, and for every function 
$\Phi\in\mathcal{O}_c(U,\mathcal{B}(C(\X,\M)))$, we may 
construct the operator $\Phi(T)\in\mathcal{B}^{\rm r}(C(\X,\H))$, using
Theorem \ref{H_afc}.
\end{Exa}

\section{Quaternionic Joint Spectrum of Paires}
\label{QJSP}

 In many applications, it is more convenient to work
with matrix quaternions rather than with abstract quaternions.
Specifically,  one considers the injective unital algebra morphism

$$
\H\ni x_1+y_1{\bf j}+x_2{\bf k}+y_2{\bf l}\mapsto
\left(\begin{array}{cc} x_1+iy_1 & x_2+iy_2 \\ -x_2+iy_2 & x_1-iy_1 \end{array}\right)\in\M_2,
$$
with  $x_1,y_1,x_2,y_2\in\R,$ where $\M_2$ is the complex algebra of $2\times2$-matrix, whose 
image, denoted by $\H_2$ is the real algebra of matrix quaternions. The elements of $\H_2$ can be also written as
matrices of the form

$$
Q(\z)=\left(\begin{array}{cc} z_1 & z_2 \\ -\bar{z}_2 & \bar{z_1} \end{array}\right),\,\,\z=(z_1,z_2)\in\C^2.
$$

A strong connection between the spectral theory of pairs of 
commuting operators in a complex Hilbert space  and the algebra of quaternions has been firstly noticed in \cite{Vas1}. Another
connection will be presented in this section. 

If $\V$ is an arbitrary vector space, we denote by $\V^2$ the
Cartesian product $\V\times\V$.

 Let $\mathcal{V}$ be a real Banach space, and let 
${\bf T}=(T_1,T_2)\in\mathcal{B(V)}^2$ be a pair of commuting
operators. The extended pair
${\bf T}_\C=(T_{1\C},T_{2\C})\in\mathcal{B(V_\C)}^2$ also consists of commuting operators. For simplicity, we set 
$$
Q({\bf T}_\C):=\left(\begin{array}{cc} T_{1\C} & T_{2\C} \\ -T_{2\C} & T_{1\C}
\end{array}\right)
$$
which acts on the complex Banach space $\mathcal{V}_\C^2$.

We now define the quaternionic resolvent set and spectrum  for the case of a pair of operators, inspired by the previous discussion
concerning a single operator.

% Def 2

\begin{Def}\label{Q-jspectrum}\rm Let $\mathcal{V}$ be a real Banach space. For a given pair ${\bf T}=(T_1,T_2)\in\mathcal{B(V)}^2$ of commuting operators, the set of those 
$Q(\z)\in\H_2,\,\z=(z_1,z_2)\in\C^2$, such that the operator
$$T_1^2+T_2^2-2\Re{z_1}T_1-2\Re{z_2}T_2+\vert z_1\vert^2+
\vert z_2\vert^2
$$
is invertible in $\mathcal{B(V)}$
is said to be the {\it quaternionic joint resolvent} (or simply the $Q$-{\it joint resolvent}) of ${\bf T}$, and is denoted by
$\rho_\H({\bf T})$. 

The complement $\sigma_\H({\bf T})=\H_2\setminus\rho_\H({\bf T})$ is called 
the {\it quaternionic joint spectrum} (or simply the $Q$-{\it joint spectrum}) of ${\bf T}$.
\end{Def}

For every pair ${\bf T}_\C=(T_{1\C},T_{2\C})\in\mathcal{B(V_\C)}^2$ we put
${\bf T}_\C^c=(T_{1\C},-T_{2\C})\in\mathcal{B(V_\C)}^2$, and for every pair 
$\z=(z_1,z_2)\in\C^2$ we put $\z^c=(\bar{z}_1,-z_2)\in\C^2$

% Lemma 3

\begin{Lem}\label{Q-jsp} A matrix quaternion $Q(\z)$\,$(\z\in\C^2)$ is in the set $\rho_\H({\bf T})$
if and only if the operators $Q({\bf T}_\C)-Q(\z),\,Q({\bf T}_\C^c)-Q(\z^c)$ are invertible in $\mathcal{B}(\mathcal{V}_\C^2)$.
\end{Lem}

{\it Proof}\, The assertion follows from the equalities
$$
\left(\begin{array}{cc} T_{1\C}-z_1 & T_{2\C}-z_2 \\
-T_{2\C}+\bar{z}_2 & T_{1\C}-\bar{z}_1\end{array}\right)
\left(\begin{array}{cc} T_{1\C}-\bar{z}_1 & -T_{2\C}+z_2 \\
T_{2\C}-\bar{z}_2 & T_{1\C}-z_1\end{array}\right)=
$$
$$
\left(\begin{array}{cc} T_{1\C}-\bar{z}_1 & -T_{2\C}+z_2 \\
T_{2\C}-\bar{z}_2 & T_{1\C}-z_1\end{array}\right)
\left(\begin{array}{cc} T_{1\C}-z_1 & T_{2\C}-z_2 \\
-T_{2\C}+\bar{z}_2 & T_{1\C}-\bar{z}_1\end{array}\right)=
$$
$$
[(T_{1\C}-z_1)(T_{1\C}-\bar{z}_1)+ 
(T_{2\C}-z_2)(T_{2\C}-\bar{z}_2)]{\bf I}.
$$
for all $\z=(z_1,z_2)\in\C^2$, where $\bf I$ is the identity.
Consequently, the operators $Q({\bf T}_\C)-Q(\z),\,Q({\bf T}_\C^c)-Q(\z^c)$ are invertible in
$\mathcal{B}({\mathcal V}_\C^2)$ if and only if the operator 
$(T_{1\C}-z_1)(T_{1\C}-\bar{z}_1)+ 
(T_{2\C}-z_2)(T_{2\C}-\bar{z}_2)$
is invertible in  $\mathcal{B}(\mathcal{V}_\C)$.  Because we have
$$
 T_{1\C}^2+T_{2\C}^2-2\Re{z_1}T_{1\C}-2\Re{z_2}T_{2\C}+\vert z_1\vert^2+\vert z_2\vert^2=
$$
$$ 
[T_1^2+T_1^2-2\Re{z_1}T_1-2\Re{z_2}T_2+\vert z_1\vert^2+
\vert z_2\vert^2]_\C,
$$
the operators $Q({\bf T}_\C)-Q(\z),\,Q({\bf T}_\C^c)-Q(\z^c)$ are invertible in
$\mathcal{B}({\mathcal V}_\C^2)$ if and only if the operator 
$T_1^2+T_1^2-2\Re{z_1}T_1-2\Re{z_2}T_2+\vert z_1\vert^2+
\vert z_2\vert^2$
is invertible in $\mathcal{B(V)}$.
\medskip

Lemma \ref{Q-jsp} shows that we have   the
property  $Q(\z)\in\sigma_\H({\bf T})$ if and only if $Q(z^c)\in\sigma_\H({\bf T}^c)$. Putting 
$$
\sigma_{\C^2}({\bf T}):=\{\z\in\C^2;Q(\z)\in\sigma_\H({\bf T})\},
$$ 
the set $\sigma_{\C^2}({\bf T})$ has a similar property, specifically
$\bf z\in\sigma_{\C^2}({\bf T})$ if and only if $\bf z^c\in\sigma_{\C^2}({\bf T}^c)$. As in the quaternionic case, the set $\sigma_{\C^2}({\bf T})$
looks like a ''complex border`` of the set  $\sigma_\H({\bf T})$.

% Rem 9

\begin{Rem}\rm For the extended pair 
${\bf T}_\C=(T_{1\C},T_{2\C})\in {B(V_\C)}^2$ of the commuting pair ${\bf T}=(T_1,T_2)\in \mathcal{B(V)}$ there is an 
interesting connexion with the {\it joint spectral theory} of 
J. L. Taylor (see \cite{Tay,Tay2}; see also \cite{Vas3}). Namely, if the operator
$T_{1\C}^2+T_{2\C}^2-2\Re{z_1}T_{1\C}-2\Re{z_2}T_{2\C}+\vert z_1\vert^2+\vert z_2\vert^2$ is invertible, then the point
$\z=(z_1,z_2)$ belongs to the joint resolvent of ${\bf T}_\C$. 
Indeed, setting 
$$
R_j({\bf T}_\C,\z)=(T_{j\C}-\bar{z}_j)(T_{1\C}^2+T_{2\C}^2-2\Re{z_1}T_{1\C}-2\Re{z_2}T_{2\C}+\vert z_1\vert^2+\vert z_2\vert^2)^{-1},
$$ $q=Q(\z)$
for $j=1,2$, we clearly have
$$
(T_{1\C}-z_1)R_1({\bf T}_\C,\z)+(T_{2\C}-z_2)R_2({\bf T}_\C,\z)
={\bf I},
$$ 
which, according to \cite{Tay}, implies that $\z$ is in the 
joint resolvent of ${\bf T}_\C$. A similar argument shows that,
in this case the point  $\z^c$ belongs to the joint resolvent of ${\bf T}_\C^c$. In addition, if $\sigma(T_\C)$
designates the Taylor spectrum of $T_\C$, we have the 
inclusion $\sigma(T_\C)\subset\sigma_{\C^2}({\bf T})$.                                                                                                                                                                                                                                                                                                                                                                                                                                                     
In particular, for every complex-valued function $f$ analytic in a  neighborhood of $\sigma_{\C^2}({\bf T})$, the operator 
$f(\bf T_\C)$ can be computed via Taylor's analytic functional
calculus. In fact, we have a Martinelli type formula for the
analytic functional calculus:
\end{Rem}

\begin{Thm} Let $\mathcal{V}$ be a real Banach space, let ${\bf T}=(T_1,T_2)\in \mathcal{B(V)}^2$ be a pair of commuting operators, let $U\subset\C^2$ be an open set, let  
$D\subset U$ be a bounded domain
containing $\sigma_{\C^2}({\bf T})$, with piecewise-smooth boundary $\Sigma$, and let $f\in\mathcal{O}(U)$. Then we have
$$
f({\bf T}_\C)=\frac{1}{(2\pi i)^2}\int_\Sigma f(\z))L({\bf z,T_\C})^{-2}(\bar{z}_1-T_{1\C})d\bar{z}_2-(\bar{z}_2-T_{2\C})
d\bar{z}_1]dz_1 dz_2,
$$
where 
$$
L({\bf z,T_\C})=T_{1\C}^2+T_{2\C}^2-2\Re{z_1}T_{1\C}-2\Re{z_2}T_{2\C}+\vert z_1\vert^2+\vert z_2\vert^2.
$$
\end{Thm}

{\it Proof.}\, Theorem III.9.9 from \cite{Vas3} implies that the 
map $\mathcal{O}(U)\ni f\mapsto f({\bf T}_\C)\in\mathcal{B(V_\C)}$, defined in terms of Taylor's analytic functional calculus, is unital, 
linear, multiplicative, and ordinary complex polynomials in $\z$
are transformed into polynomials in ${\bf T}_\C$ by simple substitution, where $\mathcal{O}(U)$ is the algebra of all analytic functions in the open set $U\subset\C^2$, provided 
$U\supset\sigma({\bf T}_\C)$.

The only thing to prove is that, when $U\supset\sigma_{\C^2}({\bf T})$, Taylor's functional calculus is given by the stated (canonical) formula. In order to do that, we use an argument  from the proof of Theorem 
III.8.1 in \cite{Vas3}, to make explicit the integral III(9.2) from \cite{Vas3} (see also \cite{Lev}).

We consider the exterior algebra
$$
\Lambda[e_1,e_2,\bar{\xi_1},\bar{\xi_2},\mathcal{O}(U)\otimes\mathcal{V}_\C]=
\Lambda[e_1,e_2,\bar{\xi_1},\bar{\xi_2}]\otimes\mathcal{O}(U)\otimes\mathcal{V}_\C,
$$
where the indeterminates $e_1,e_2$ are to be associated with the pair ${\bf T}_\C$, we put $\bar{\xi_j}=d\bar{z}_j,\,j=1,2$, and consider the operators
$\delta=(z_1-T_{1\C})\otimes e_1+(z_2-T_{2\C})\otimes e_2,\,\bar{\partial} = 
(\partial/\partial\bar{z_1})\otimes\bar{\xi_1}+(\partial/\partial\bar{z_2})\otimes\bar{\xi_2}$, acting naturally on this 
exterior algebra, via the calculus with exterior forms. 

To simplify the computation, we omit the symbol $\otimes$, and the exterior product will be denoted simply par juxtaposition. 

We fix the exterior form $\eta=\eta_2=fye_1e_2$ for some 
$f\in\mathcal{O}(U)$ and  $y\in\mathcal{X}_\C$, which clearly satisfy the equation 
$(\delta+\bar{\partial})\eta=0$, and look for a solution
$\theta$ of the equation $(\delta+\bar{\partial})\theta=\eta$. 
We write $\theta=\theta_0+\theta_1$, where $\theta_0,\theta_1$
are of degree $0$ and $1$ in $e_1,e_2$, respectively. Then the equation $(\delta+\bar{\partial})\theta=\eta$ can be written 
under the form $\delta\theta_1=\eta,\,\delta\theta_0=-\bar{\partial}\theta_1$, and $\bar{\partial}\theta_0=0$. 
Note that
$$
\theta_1=fL({\bf z,T_\C})^{-1}[(\bar{z}_1-T_{1\C})ye_2-
(\bar{z}_2-T_{2\C})]ye_1
$$
is visibly a solution of the equation $\delta\theta_1=\eta$. 
Further, we have
$$
\bar{\partial}\theta_1=fL({\bf z,T_\C})^{-2}
[(z_1-T_{1\C})(\bar{z}_2-T_{2\C})y\bar{\xi}_1e_1-
(z_1-T_{1\C})(\bar{z}_1-T_{1\C})y\bar{\xi}_2e_1+
$$
$$
(z_2-T_{2\C})(\bar{z}_2-T_{2\C})y\bar{\xi}_1e_2-
(z_2-T_{2\C})(\bar{z}_1-T_{1\C})y\bar{\xi}_2e_2]=
$$
$$
\delta[fL({\bf z,T_\C})^{-2}(\bar{z}_1-T_{1\C})y\bar{\xi}_2-
fL({\bf z,T_\C})^{-2}(\bar{z}_2-T_{2\C})y\bar{\xi}_1],
$$
so we may define
$$
\theta_0=-fL({\bf z,T_\C})^{-2}(\bar{z}_1-T_{1\C})y\bar{\xi}_2+
fL({\bf z,T_\C})^{-2}(\bar{z}_2-T_{2\C})y\bar{\xi}_1.
$$
Formula III(8.5) from \cite{Vas3} shows that 
$$
f({\bf T}_\C)y=-\frac{1}{(2\pi i)^2}\int_U\bar{\partial}(\phi\theta_0)dz_1 dz_2=
$$
$$
\frac{1}{(2\pi i)^2}\int_\Sigma f(\z))L({\bf z,T_\C})^{-2}[(\bar{z}_1-T_{1\C})yd\bar{z}_2-(\bar{z}_2-T_{2\C})y
d\bar{z}_1]dz_1 dz_2,
$$
for all $y\in\mathcal{X}_\C$, via Stokes's formula, where $\phi$ is a smooth function such 
that $\phi=0$ in a neighborhood of $\sigma_{\C^2}({\bf T})$,
$\phi=1$ on $\Sigma$ and the support of $1-\phi$ is compact.

% Rem 10

\begin{Rem}\rm (1) We may extend the previous functional calculus to $\mathcal{B(V}_\C)$-valued analytic functions, setting, for such a function $F$ and with the notation from above,
$$
F({\bf T}_\C)=\frac{1}{(2\pi i)^2}\int_\Sigma F(\z))L({\bf z,T_\C})^{-2}(\bar{z}_1-T_{1\C})d\bar{z}_2-(\bar{z}_2-T_{2\C})
d\bar{z}_1]dz_1 dz_2.
$$
In particular, if $F(\z)=\sum_{j,k\ge0}A_{jk\C}z_1^jz_2^k$,
with $A_{j,k}\in\mathcal{B(V)}$, where the series is convergent in neighborhood of $\sigma_{\C^2}({\bf T})$, we obtain 

$$F({\bf T}):=F({\bf T}_\C)\vert\mathcal{V}=\sum_{j,k\ge0}A_{jk}T_1^jT_2^k\in\mathcal{B(V)}.$$
 
(2) The connexion of the spectral theory of pairs with the algebra of quaternions is even stronger in the case of complex Hilbert
spaces. Specifically, if $\mathcal{H}$ is a complex Hilbert space and ${\bf V}=(V_1,V_2)$ is a commuting pair of bounded
linear operators on $\mathcal{H}$, a point $\z=(z_1,z_2)\in\C^2$
is in the joint resolvent of ${\bf V}$ if and only if the 
operator $Q({\bf V})-Q(\z)$ is invertible in $\mathcal{H}^2$,
where 
$$
Q({\bf V})=\left(\begin{array}{cc} V_1 & V_2 \\ -V_2^* & V_1^*
\end{array}\right).
$$
(see \cite{Vas1} for details). In this case, there is also a Martinelli type formula which can be used to construct the associated analytic functional calculus (see \cite{Vas2},\cite{Vas3}). An approach to such a construction in Banach spaces, by using a so-called splitting joint spectrum, can be 
found in \cite{MuKo}.  
\end{Rem}

\end{document}